\documentclass[12pt,letterpaper,final, leqno, twoside]{article}
\usepackage{amsmath}

\usepackage{amsthm}
\usepackage{amssymb}
\usepackage{amscd}
\usepackage{graphicx}
\usepackage{afterpage}

\begin{document}

\title    {\ \textbf{Fixed point properties of nilpotent and solvable Lie group actions on Hadamard manifolds} }
\author  {\small   Mehrzad Monzavi}
\date{\today}
\maketitle {\footnotesize
 \vskip
0.2 true cm \begin{abstract}
We will prove the following theorems. The first theorem posits the existence of a fixed point for the actions of nilpotent Lie groups on nonpositively curved compact manifolds. The second theorem states that actions of solvable Lie groups on nonpositively curved compact manifolds have either a fixed point or a $2${\textendash}periodic point.\\


 \vskip 0.2
true cm \vskip 0.2 true cm {\it Key words:} Fixed Point Property, Lie Group Action, Hadamard Manifold.
\end{abstract}}
 \vskip
0.2 true cm

\begin{flushleft}
  1. {\bf \small Introduction}
\end{flushleft}

Let $X$ be a topological space, $G$ a group, and $\phi : G \times X \rightarrow X$ a continuous action of $G$ on $X$. The set $Gx=\{gx : g\in G\}$ is called the orbit of $x$ under the action of $G$. For $x\in X$, If $Gx$ consists of $n$ points, $x$ is called an $n${\textendash}periodic point, . A $1${\textendash}periodic point $x$ is called a fixed point. In other words a point $x$ is called a fixed point if $gx = x$ for all $g \in G$.The following question is of high interest:\\

Under which conditions on $G$ and $X$, will any  $\phi : G \times X \rightarrow X$ have fixed points or periodic points?\\

A Continuum (compact, connected, metric space) $X$ is said to be $1${\textendash}arcwise connected if given any two points $x, y \in X$, $x \neq y$, there is a unique arc in $X$ whose endpoints are $x$ and $y$. In other words, $X$ is arcwise connected and contains no simple closed curves.\\

It is proved in \cite{R2} that If $X$ is a $1${\textendash}arcwise connected continuum, and $G$ is the discrete cyclic group $\mathbb{Z}$, then $\phi$ has a fixed point.In other words, every homeomorphism $h: X \rightarrow X$ leaves some point fixed.\\ 

It is proved in \cite{R4} that If $X$ is a $1${\textendash}arcwise connected continuum, and $G$ is a nilpotent group then $\phi$ has a fixed point.\\

It is proved in \cite{R5} that If $X$ is a $1${\textendash}arcwise connected continuum, and $G$ is a solvable group then $\phi$ has a fixed point or a $2${\textendash}Periodic point.\\

A Hadamard manifold $H$ is a complete, simply connected Riemannian manifold of dimension $n\geq2$ having sectional curvature $K \leq 0$.\\ 

A Hadamard manifold $H$ can be compactified by adding the visual boundary $\partial{H}$. (\cite{R1})\\

This paper will generalize these results to Lie group actions on Hadamard manifolds.\\

We will prove the following theorems:\\

{\textbf{Theorem 1.1.}} {\it Suppose $H$ is a compact Hadamard manifold and $N$ is a nilpotent Lie group. If $\phi:N\times H \rightarrow H$ is a continous action, then $\phi$ has a fixed point.\\}

{\textbf{Theorem 1.2.}} {\it Suppose $H$ is a compact Hadamard manifold and $S$ is a solvable Lie group. If $\phi:S\times H \rightarrow H$ is a continuous action, then $\phi$ has either fixed point or a $2${\textendash}periodic point.\\}

 Proofs of these theorems use the strategy developed in \cite{R4} and \cite{R5}.\\

\newpage

{\bf 2. Preliminaries}\\\\

We will use the following definitions.\\

Let $G$ be a Lie group with and $[a,b] = a^{-1}b^{-1}ab$ the commutator.\\

(1) Consider the lower central sequence $G = G_0 \triangleright G_1 \triangleright G_2 \ldots$ where $G_{k+1} = [G,G_k]$ and $G_k/G_{k+1}$ is Abelian.\\

A group $G$ is called nilpotent if there exists $n$ such that $G_{n}= \{e\}$\\

(2)A Lie group is said to be nilpotent if it is nilpotent as an abstract group.\\

(3) Consider the normal subgroup series $G = G^0 \triangleright G^1 \triangleright G^2 \ldots$ where $G_{k+1} = [G^k,G^k]$ and $G_k/G_{k+1}$ is Abelian.\\~\\

A group $G$ is called solvable if there exists $n$ such that $G^{n}= \{e\}$\\

(4) A Lie group is said to be solvable if it is solvable as an abstract group.\\

(5) The symbol $[x,y]$ in a Hadamard manifold $H$ denotes the unique geodesic connecting $x$ to $y$.\\

(6) $P_n(G, X)$ denotes all n{\textendash}periodic points of $G$ and $P_{\leq n}(G,X)$ denotes all $n${\textendash}periodic points with periods less than or equal to $n$. In other words, $P_{\leq n}(G,X) = \bigcup_{i=1}^{n} P_n(G,X)$.\\\\

{\bf \large {2. Proof of Theorem 1.1.}}\\\\

{\textbf{Lemma 2.1.}} {\it Let $H$ be a compact Hadamard manifold and $\phi:N\times H \rightarrow H$ be an action of Nilpotent Lie group $N$ on $H$. Suppose $K$ is a normal subgroup of $N$ and ${N}/{K}$ is a cyclic group. If $P_{1}(K,H) \neq \emptyset$ then $P_{1}(N,H) \neq \emptyset$.\\}

\textbf{Proof of Lemma 2.1.}\\

Let ${N}/{K}=<\bar{g}_0>$ for some $g_0 \in N$. $g_0(P_1(K, H)) = P_1(K, H)$\\

Let $\mathcal{A}$ = \{ $[x,g_0(x)]$ : $x \in P_1(K,H) \}$. Define a partial order $\prec$ in $\mathcal{A}$ as follows:\\

$[x,g_0(x)]$ $\prec$ $[x',g_0(x')]$ iff $[x,g_0(x)]$ $\supseteq$ $[x',g_0(x')]$.\\

Since every totally ordered subset of $\mathcal{A}$ has an upper bound which is the intersection of all its elements, by Zorn's lemma, there exists a maximal geodesic $[y_0,g_0(y_0)] \in \mathcal{A}$.\\


We will have the following cases:\\

\textbf{Case i.} $g_0(y_0)=y_0$. In this case, $y_0$ is a fixed point of $N$.\\

\textbf{Case ii.} ${g_0}^2(y_0) = y_0$. In this case, We will show that a continuous action of a nilpotent Lie group $N$ on $[y_0,g_0(y_0)]$ has a fixed point.\\

Suppose $y_0$ and $g_0(y_0)$ are not fixed points of $N$. Therefore there is a $g_1 \in N$ such that $g_1(y_0) = g_0(y_0)$ and $g_1(g_0(y_0)) = y_0$. Thus $g_1$ has a unique fixed point $x_0 \in [y_0,g_0(y_0)]$. It will be shown that $x_0$ is a fixed point of $N$.\\

Consider the following sequence;\\

$N_0=\{e\}$ and $N_{i+1} = \{ g \in N : [g,g_0] \in N_i \}$.\\ 

$\{e\} = N_0 \subset N_1 \subset \ldots$\\

Since $N$ is a Nilpotent Lie group, this sequence termintates with $N$ after $m$ stops. If $x_0$ is a fixed point of $N_i$, for $g \in N_{i+1}$, $gx_0$ is a fixed point of $g_1$ because;\\

$g_1(gx_0) = gg^{-1}{g_1}^{-1}g(g_1x_0) = g[g,g_1]x_0 = gx_0$\\

Since $x_0$ is the unique fixed point of $g_1$, $gx_0 = x_0$. Thus $x_0$ is a fixed point of $N_{i+1}$. Inductively, $x_0$ is a fixed point of $N$.\\

\textbf{Case iii.} $g_0(y_0) \neq y_0$ and ${g_0}^2(y_0) \neq y_0$\\


$[{g_0}^{n−1}(y) ,{g_0}^n(y) ] \cap [{g_0}^n(y) ,g^{n+1}(y) ]={g^n(y)}$, for all $n \in \mathbb{Z}$. (\cite{R2})\\

Consider;
\begin{center}

$ \ldots \cup [{g_0}^{-1}(y_0), y_0)] \cup [y_0,g_0(y_0)] \cup [g_0(y_0),{g_0}^2(y_0)] \cup \ldots$\\

\end{center}

By a measure theory argument from \cite{R2}, one of the following cases happens:\\

\begin{center}

$[y_0,g_0(y_0)] \cup [g_0(y_0),{g_0}^2(y_0)] \cup \ldots \subseteq [a,b]$\\~\\

$ \ldots [{g_0}^{-2}(y_0), {g_0}^{-1}(y_0)] \cup [{g_0}^{-1}(y_0), y_0)] \subseteq [a,b]$\\~\\

\end{center}

 where $[a,b]$ is a geodesic in $H$\\

In the first case, the sequence ${g_0}^n(y_0)$ converges to a fixed point. (\cite{R2})\\

Similarly, in the second case, the sequence ${g_0}^{-n}(y_0)$ converges to a fixed point. (\cite{R2})\\

{\textbf{Lemma 2.2.}} {\it Let $H$ be a compact Hadamard manifold and $\phi:N\times H \rightarrow H$ be an action of Nilpotent Lie group $N$ on $H$. Suppose $K$ is a normal subgroup of $N$ and ${N}/{K}$ is a finitely generated Abelian group. If $P_{1}(K, H) \neq \emptyset$ then $P_{1}(N,H) \neq \emptyset$.\\}

\textbf{Proof of Lemma 2.2.}\\

Since $N/K$ is a finitely generated Abelian group, $N/K \cong \mathbb{Z}_{l_1} \bigoplus \ldots \bigoplus \mathbb{Z}_{l_n}$ where  $\mathbb{Z}_{l_i} = \mathbb{Z}/<l_i>$  for $1 \leq i \leq n$. ($l_i$ may be $0$)\\

Let $\bar{K}_i < N/K$ suth that $\bar{K}_i = \mathbb{Z}_{l_1} \bigoplus \ldots \bigoplus \mathbb{Z}_{l_i}$\\

Let $\phi: N \rightarrow N/K$ be the quotient homomorphisom. \\

Let $K_i = \phi^{-1} (\bar{K}_i)$;\\

\begin{center}
$K_{i+1}/K_i \cong (K_{i+1}/K) / (K_{i}/K)$ = \\~\\

$\mathbb{Z}_{l_1} \bigoplus \ldots \bigoplus \mathbb{Z}_{l_{i+1}} /  \mathbb{Z}_{l_1} \bigoplus \ldots \bigoplus \mathbb{Z}_{l_{i}} \cong \mathbb{Z}_{l_{i+1}}$\\~\\

\end{center}

Thus $(K_{i+1}/K_i)$ is a cyclic group.\\

Consider $K \triangleright K_1 \triangleright K_2 \ldots \triangleright K_n = G$;\\

Since $P_{1}(K,H) \neq \emptyset$, by using Lemma $2.1$ repeadedly, $P_{1}(N,H) \neq \emptyset$.\\

{\textbf{Lemma 2.3.}} {\it Let $H$ be a compact Hadamard manifold and $\phi:N\times H \rightarrow H$ be a continuous action of Nilpotent Lie group $N$ on $H$. Suppose $K$ is a normal subgroup of $N$ and ${N}/{K}$ is Abelian. If $P_{1}(K,H) \neq \emptyset$ then $P_{1}(N,H) \neq \emptyset$.\\}

\textbf{Proof of Lemma 2.3.}\\

Consider the group generated by $K$ and $T$ where $T$ is a finite subset of $N$.

$<K,T>/K$ is a finitely generated Abelian group.\\

By Lemma $2.2.$, $P_{1}(<K,T>) \neq \emptyset$.\\

Let $\mathcal{C} =  \{P_{1}(<K,T>) : S \subset G, T$ finite $\}$\\

Any finite intersection of $\mathcal{C}$ has a non empty intersection. In other words, $\mathcal{C}$ has the finite intersection property. \\

By compactness of $H$, $\bigcap_{c \in \mathcal{C}} c \neq \emptyset$.\\

Since $\bigcap_{c \in \mathcal{C}} c \subseteq P_1(N,H)$, $P_1(N,H) \neq \emptyset$\\


\newpage

\textbf{Proof of Theorem 1.1.}\\

Consider the lower central sequence $N = N_0 \triangleright N_1 \triangleright \ldots \triangleright N_n = \{e\}$ of $N$. $N_{i+1} = [N,N_i]$ and $N_i/N_{i+1}$ is an Abelian group, for $0 \leq i \leq n - 1$. Because $P_{1}(N_n,H)= P_{1}(\{e\})= H \neq \emptyset$, it follows inductively
from Lemma $2.3$ that $P_{1}(N,H) \neq \emptyset$.\\\\

{\bf \large {3. Proof of Theorem 1.2}}\\\\

{\textbf{Lemma 3.1.}} {\it Let $H$ be a compact Hadamard manifold and $\phi:S\times H \rightarrow H$ be a continous action of Solvable Lie group $S$ on $H$. Suppose $K$ is a normal subgroup of $S$ and ${S}/{K}$ is a cyclic  group. If $P_{\leq 2}(K,H) \neq \emptyset$ then $P_{\leq 2}(S,H) \neq \emptyset$.\\}

\textbf{Proof of Lemma 3.1.}\\

Let $\mathcal{A}$ = \{ $[x,g_0(x)]$ : $x \in P_{\leq 2}(K,H) \}$, by a similar argument to that of Lemma $2.1.$, there exists a maximal element $[x,g_0(x)] \in \mathcal{A}$.\\

If $g_0(x)=x$, $x \in P_{1}(S,H)$. Suppose ${g_0}(x) \neq x$;\\

\textbf{Case i.} $x \in P_1(K,H)$  \\

Since $[z,g_0(x)] = [x,g_0(x)] \cap [g_0(x), {g_0}^2(x)]$, $k([z,g_0(x)]) = k([x,g_0(x)]) \cap k([g_0(x), {g_0}^2(x)]) = [x,g_0(x)] \cap [g_0(x), {g_0}^2(x)] = [z,g_0(x)]$ for every $k \in K$. Therefore, $z \in P_1(K,H)$.\\


By maximality of $[x,g_0(x)]$ and $[g_0(x),{g_0}^2(z)]$ , $z \notin (x, g_0(x))$ and $z \notin (g_0(x), {g_0}^{2}(x))$; If $z \in (x, g_0(x))$, ${g_0}^{-1}(z)$ is also a fixed point by using the fact that $g_0$ maps fixed points of $K$ to fixed points of $K$ and maps $2$-periodic points of $K$ to $2$-periodic points of $K$. But $[{g_0}^{-1}(z),z] \subsetneq [x,g_0(x)]$. Also, if $z \in (g_0(x), {g_0}^{2}(x))$, then $[z,{g_0}(z)] \subsetneq [g_0(x),{g_0}^2(z)]$. \\

We will have the following subcases:\\

\textbf{Subcase i.1.} $z = g_0(x) = {g_0}^2(x)$. In this subcase, $x \in P_1(S,H)$ because $[x,{g_0}(x)]$ is invarinat under the action of $N$. So $x \in P_{\leq 2}(S,H)$.\\

\textbf{Subcase i.2.} $z = g_0(x)$. By a measure theory argument from \cite{R2}, $P_{\leq 2}(S,H) \neq \emptyset$.\\

\textbf{Case ii.} $x \in P_2(K,H)$\\

$Nx = \{x , x'\}$ therefore $Ng_0(x) = g_0(Nx) = \{g_0(x), g_0(x')\}$. Suppose $x, x', g_0(x)$ and $g_0(x')$ are pairwise different. Thus $[x, x'] \cap [g_0(x), g_0(x')] \neq \emptyset$ or $[x, g_0(x)] \cap [x', g_0(x')] \neq \emptyset$\\

By maximality of $[x,g_0(x)]$, the case that $[x',x] \cap [g_0(x'),g_0(x)] = \emptyset$ and $[x',g_0(x')] \cap [x,g_0(x)] = [a,b]$ for some $a,b \in H$ would not happen; $[x',x]$ and $[g_0(x'),g_0(x)]$ will both be $K$-invariant, thus $a$ will be the unique $K$-fixed point in $[x',x]$ which would be mapped to $b$, the unique $K$-fixed point in $[g_0(x'),g_0(x)]$, by $g_0$.\\

We will have the following subcases:\\

\textbf{Subcase ii.1.} $[x',x] \cap [g_0(x'),g_0(x)] = [a,b]$ for some $a,b \in H$ and $[x',g_0(x')] \cap [x,g_0(x)] =\emptyset$. Since $[a, b]$ is $K-$invariant, $a, b \in P_{\leq 􏰶2}(K, H)$\\

By maximality of $[x,g_0(x)]$, $b \notin (x, g_0(x))$. Otherwise, $g_0(b) \in [b, g_0(x')]$. Thus $g^{-1}(b) \in [x,b]$. Therefore $[g^{-1}(b),b]􏰭 \subsetneq [x, g_0(x)]$.\\


Moreover, the case that $b = x$ and  $a \in (x', g_0(x'))$ would not also happen. Otherwise, $b = x \in P_1(K,H)$. So $a = x'$ or $a = g_0(x')$ while $a \in (x', g_0(x'))$ which yeilds to $P_{\leq 2}(S,H) \neq \emptyset$ by a measure theory argument from \cite{R2}.\\

The case that $b = g_0(x)$ similarly yeilds to $P_{\leq 2}(S,H) \neq \emptyset$.\\

\textbf{Subcase ii.2.}  $[x',x] \cap [g_0(x'),g_0(x)] = [x',g_0(x')] \cap [x,g_0(x)] =[a,b]$ for some $a,b \in H$.\\

By maximality of $[x,g_0(x)]$, $b \notin (x, g_0(x'))$. Otherwise $g_0(b) \in [b, g_0(x')]$. Thus ${g_0}^{-1}(b) \in [x,b]$. Then $[{g_0}^{-1}(b),b] \subseteq [x, g_0(x)]$. Therefore $b = g_0(x)=a$. Thus $g_0(x) \in (x, g_0(x'))$. So $g_0(x) \in P_1(K,H)$.\\ 

Furthermore, Since $x \in P_2(K,H)$, the case that $b = x$ and  $a \in (g_0(x), x')$ would not happen. So if $b = x$. Then $a = g_0(x)$ or $a =x'$.\\

The case that $b = g_0(x')$ similarly yeilds to $P_{\leq 2}(S,H) \neq \emptyset$.\\


In the case of $[x',x] \cap [g_0(x'),g_0(x)] = [x',g_0(x')] \cap [x,g_0(x)] =[a,x], a = g_0(x)$. Since $[x,g_0(x)]$ and $[g_0(x),{g_0}^{2}(x)]$ are maximal, $[x,g_0(x)] = [g_0(x),{g_0}^{2}(x)]$. Since $[x,g_0(x)]$ is $K$-invariant, $[x,{g_0}(x)]$ is also $S$-invarinat. Thus $x \in P_{\leq 2}(S,H)$.\\

{\textbf{Lemma 3.2.}} {\it Let $H$ be a compact Hadamard manifold and $\phi:S \times H \rightarrow H$ be a continuous action of Solvable Lie group $S$ on $H$. Suppose $K$ is a normal subgroup of $S$ and ${S}/{K}$ is a finitely generated Abelian group. If $P_{\leq 2}(K,H) \neq \emptyset$ then $P_{\leq 2}(S,H) \neq \emptyset$.\\}

Lemma $3.2.$ can be proved similarly to the proof of Lemma $2.2.$\\











{\textbf{Lemma 3.3.}} {\it Let $H$ be a compact Hadamard manifold and $\phi:S\times H \rightarrow H$ be a continuous action of Solvable Lie group $S$ on $H$. Suppose $K$ is a normal subgroup of $S$ and ${S}/{K}$ is Abelian. If $P_{\leq 2}(K,H) \neq \emptyset$ then $P_{\leq 2}(S,H) \neq \emptyset$.\\}


Lemma $3.3.$ can be proved similarly to the proof of Lemma $2.2.$\\~\\








\textbf{Proof of Theorem 1.2.}\\

Consider the normal subgroup series $S = S^0 \triangleright S^1 \triangleright \ldots \triangleright S^n = \{e\}$ of $S$. $S^{i+1} = [S^i,S^i]$ and $S^i/S^{i+1}$ is an Abelian group, for $0 \leq i \leq n - 1$. Because $P_{\leq 2}(S^n,H)= P_{2}(\{e\}) \neq \emptyset$, it follows inductively
from Lemma $3$ that $P_{\leq 2}(S,H) \neq \emptyset$.

 \newpage

  {\footnotesize

  \begin{center}


\end{center}

\bigskip

\begin{flushleft}

\footnotesize {{\it
 \textsc{Department of mathematics} \ \ \ \ \ \ \ \\

\textsc{University of Texas at Arlington}\\
 \textsc{Arlington, Texas, 76013, USA}\ \ \ \ \ \ \ \ \ \ \ \ \ \ \ \ \\
  \textsc{mehrzad.monzavi@mavs.uta.edu} \ \ \ \ \\

    }}

\end{flushleft}
}
\end{document}